\documentclass[12pt]{amsart}
\usepackage[a4paper,left=3cm,right=3cm,top=3cm,bottom=3cm]{geometry}
\usepackage{hyperref}
\usepackage{amssymb}
\usepackage{amsthm,amsmath}
\usepackage{enumerate}
\usepackage{graphicx}
\usepackage[OT4,T1]{fontenc}
\usepackage[utf8]{inputenc}

\usepackage{graphicx}

\usepackage{amsmath}
\usepackage{amssymb}

\usepackage{lastpage} 

\usepackage{etex} 

\usepackage{color}
\usepackage{cite} 

\RequirePackage[numbers]{natbib}

\newtheorem{thm}{Theorem}
\theoremstyle{remark}
\newtheorem{remark}{Remark}

\renewcommand{\le}{\leqslant}
\renewcommand{\ge}{\geqslant}

\begin{document}
\title[A model for random fire induced tree-grass coexistence in savannas]{A model for random fire induced tree-grass coexistence in savannas}

\author[P. Klimasara]{Paweł Klimasara}
\address{Department of Biomathematics, Institute of Mathematics, University of Silesia in Katowice, Bankowa 14, 40-007 Katowice}
\email{p.klimasara@gmail.com}

\author[M. Tyran-Kamińska]{Marta Tyran-Kamińska}
\address{Institute of Mathematics, Polish Academy of Sciences, Bankowa 14, 40-007 Katowice}
\email{mtyran@us.edu.pl}

\thanks{This work was partially supported  by the Polish NCN grant  2017/27/B/ST1/00100 and  by the grant 346300 for IMPAN from the Simons Foundation and the matching 2015-2019 Polish MNiSW fund.}


\begin{abstract}
Tree-grass coexistence in savanna ecosystems depends strongly on environmental disturbances out of which crucial is fire. Most modeling attempts in the literature lack stochastic approach to fire occurrences which is essential to reflect their unpredictability. Existing models that actually include stochasticity of fire are usually analyzed only numerically. We introduce new minimalistic model of tree-grass coexistence where fires occur according to stochastic process. We use the tools of linear semigroup theory to provide more careful mathematical analysis of the model. Essentially we show that there exists a unique stationary distribution of tree and grass biomasses.
\end{abstract}
\maketitle

\section{Introduction}
Savanna covers around 20\% of the Earth's land surface. It is a mixed woodland-grassland ecosystem with canopy open enough to support the existence of continuous herbaceous layer dominated by grass. In order to find the explanation of such tree-grass codominance many theoretical models were introduced. Beside interspecies competition (e.g. \cite{eagleson1985water}), this coexistence is believed to have been driven by various environmental disturbances, primarily rainfall (e.g. \cite{ursino2014eco}, \cite{synodinos2018impact}), grazing and browsing (e.g. \cite{bodini2016vegetation}), and fire \cite{sankaran2008woody}. Some models consider additional factors like competition of tree seedlings with grass \cite{baudena2010idealized} or varying flammability of trees \cite{beckage2009vegetation}. From the mathematical point of view, models containing many different factors lack stochasticity and differ in methodology (see e.g. loop analysis for graphs in \cite{bodini2016vegetation} or models based on impulsive differential equations \cite{tamen2016tree}, \cite{yatat2017impulsive}).

Realistically, the appearance of fire is stochastic and its frequency can vary significantly \cite{archibald2009limits}. Usually studies with stochastic fire focus on numerical analysis (see e.g. \cite{d2006probabilistic}, \cite{baudena2010idealized}, \cite{beckage2011grass}, \cite{de2014tree}, and \cite{synodinos2018impact}). We introduce a simple model where fire occurrences are stochastic and study it in terms of linear semigroup theory. We find that biomasses of grass and trees have a unique stationary distribution and hence this simple model can describe stable savannas driven by stochasticity of fires.

\section{Model description}
Our model is based on a simplified version of the system of differential equations given in \cite{beckage2011grass}, but instead of putting fire disturbances inside these equations we introduce appropriate stochastic process separately. Similarly to cited authors we consider only amounts of tree and grass biomasses, fires are events discrete in time, and the strength of grass-fire feedback depends on biomass of grass.

\begin{figure}[h!]
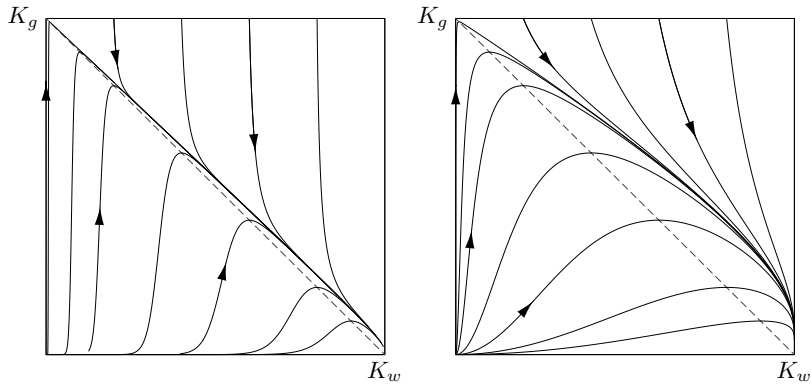

\centering
 \includegraphics[height=5cm]{mphase-1}
 \includegraphics[height=5cm]{mphase-2}
 \caption{Phase portraits for \eqref{no-fire} with parameter values $r_w=0.08$, $r_g=1.5$ (left-hand panel) and $r_w=0.25$, $r_g=0.5$ (right-hand panel), $K_w=K_g=1$}
 \label{fig-phase}
  \end{figure}

In the absence of fires we represent the dynamics of tree biomass $W$ and grass biomass $G$ (both in $\frac{\text{mass}}{\text{area}}$ units)  according to the  competition model
\begin{equation}\label{no-fire}
\left\{\begin{array}{lll}
W'(t)=r_w W(t)\left(1-\frac{W(t)}{K_w}\right),\\ G'(t)=r_q G(t)\left(1-\frac{G(t)}{K_g}-\frac{W(t)}{K_w}\right),
\end{array} \right.
\end{equation}
where $r_w, r_g$ are the growth rates and $K_w,K_g$ are the carrying capacities for tree and grass biomasses.
It is easily seen that \eqref{no-fire} has three stationary states: $(0,0)$, $(K_w,0)$ and $(0,K_g)$. Moreover, the point $(K_w,0)$ is locally stable, while the points $(0,0)$ and $(0,K_g)$ are unstable.
So the system of equations (\ref{no-fire}) provides a deterministic description of the change of wood and grass biomasses in time where in the long time, due to species competition, the system will end up as a woodland. The solution curves for the system~\eqref{no-fire} have the qualitative behavior as shown in Figure \ref{fig-phase}.

\begin{figure}[h!]
\includegraphics[height=5cm]{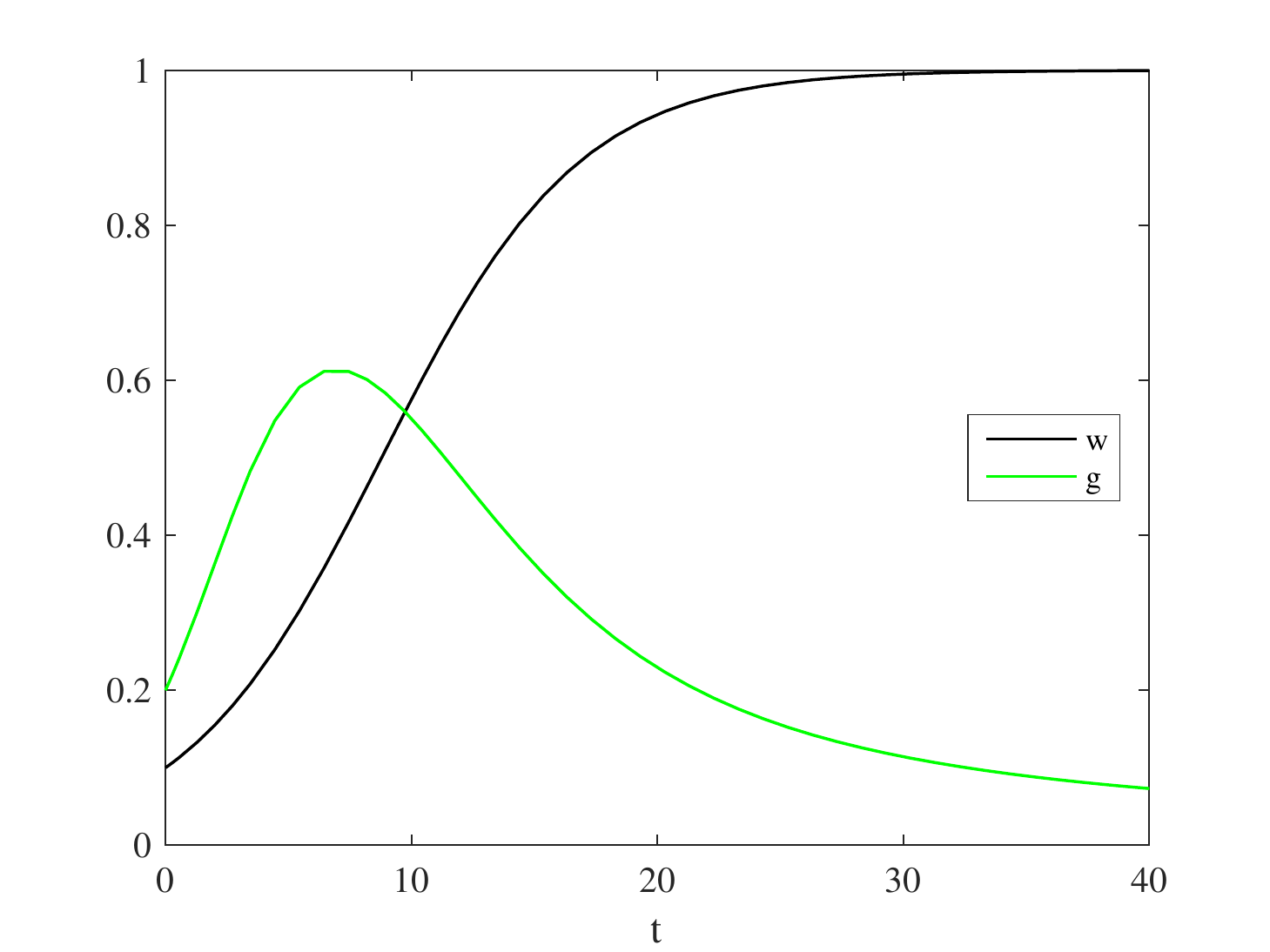}
\includegraphics[height=5cm]{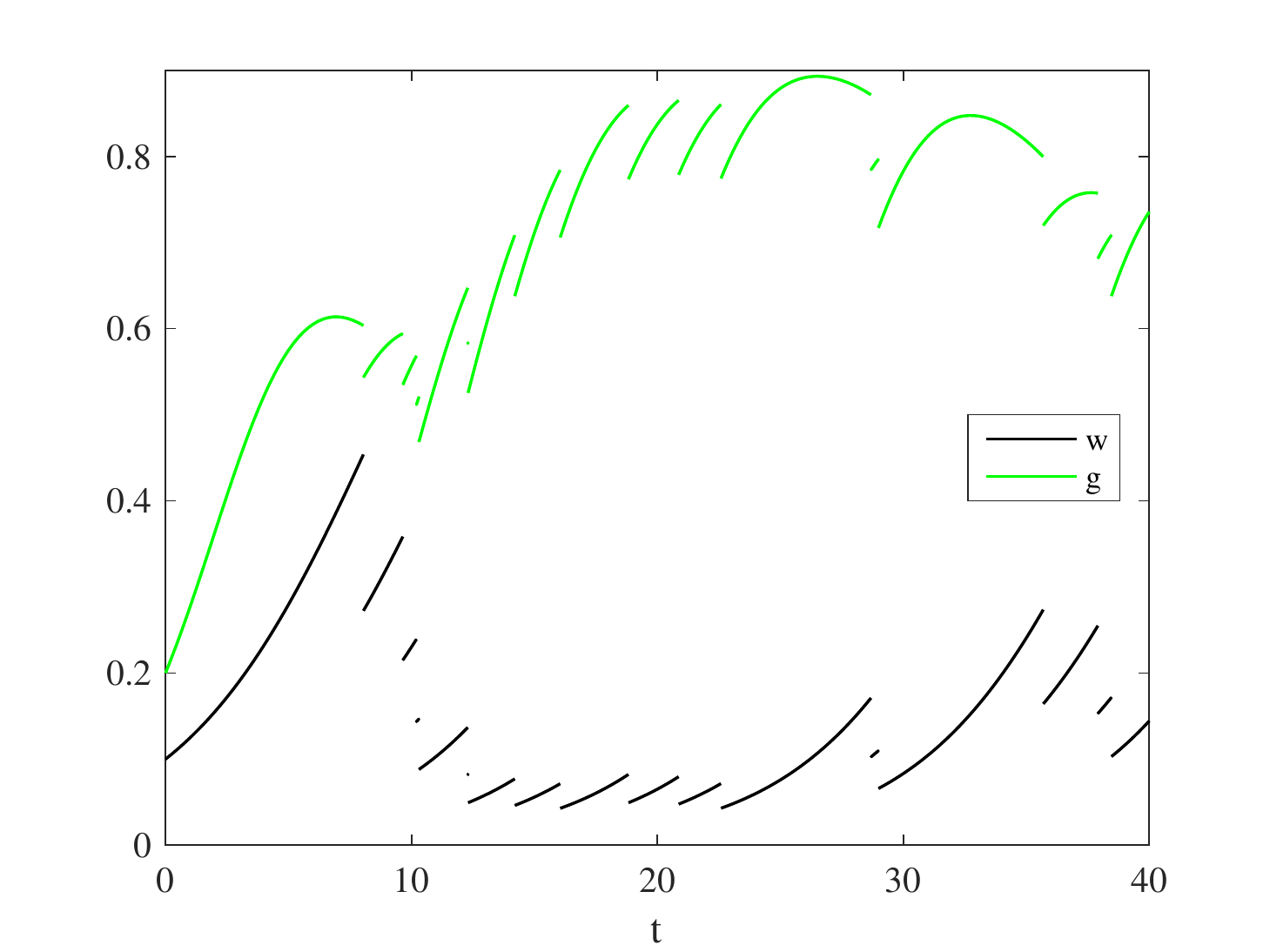}
\caption{Graphs of system \eqref{no-fire} (left-hand panel) and of system \eqref{basic.eq}--\eqref{loss.eq} (right-hand panel) with parameter values $r_w=0.25$, $r_g=0.5$, $M_w=0.4$, $M_g=0.1$, $\lambda(w,g)=g$ and initial condition $w_0=0.1,\ g_0=0.2$}
\label{fire-no-fire}
\end{figure}

\begin{figure}[h!]
\centering
 \includegraphics[height=5cm]{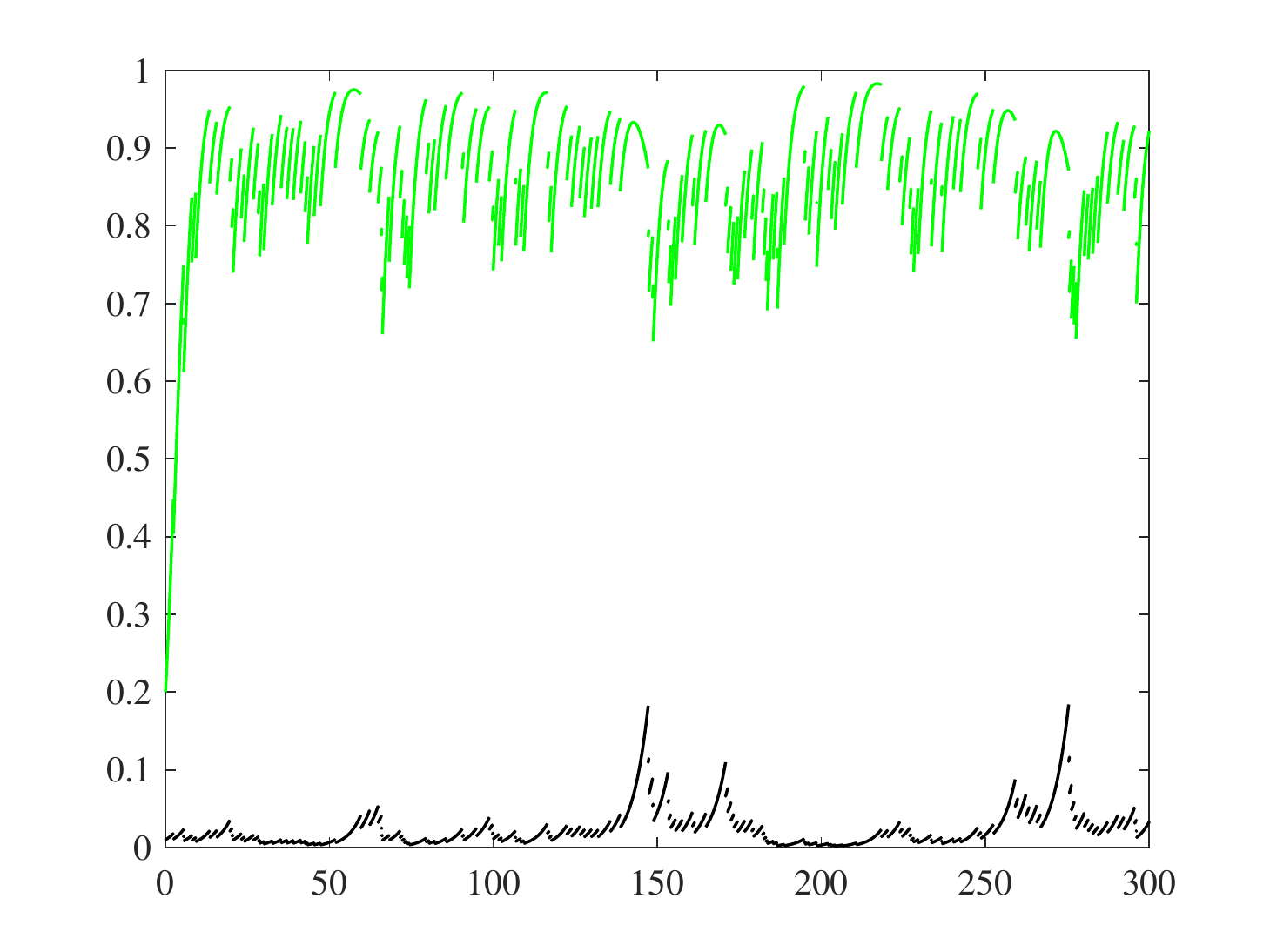}%
 \includegraphics[height=5cm]{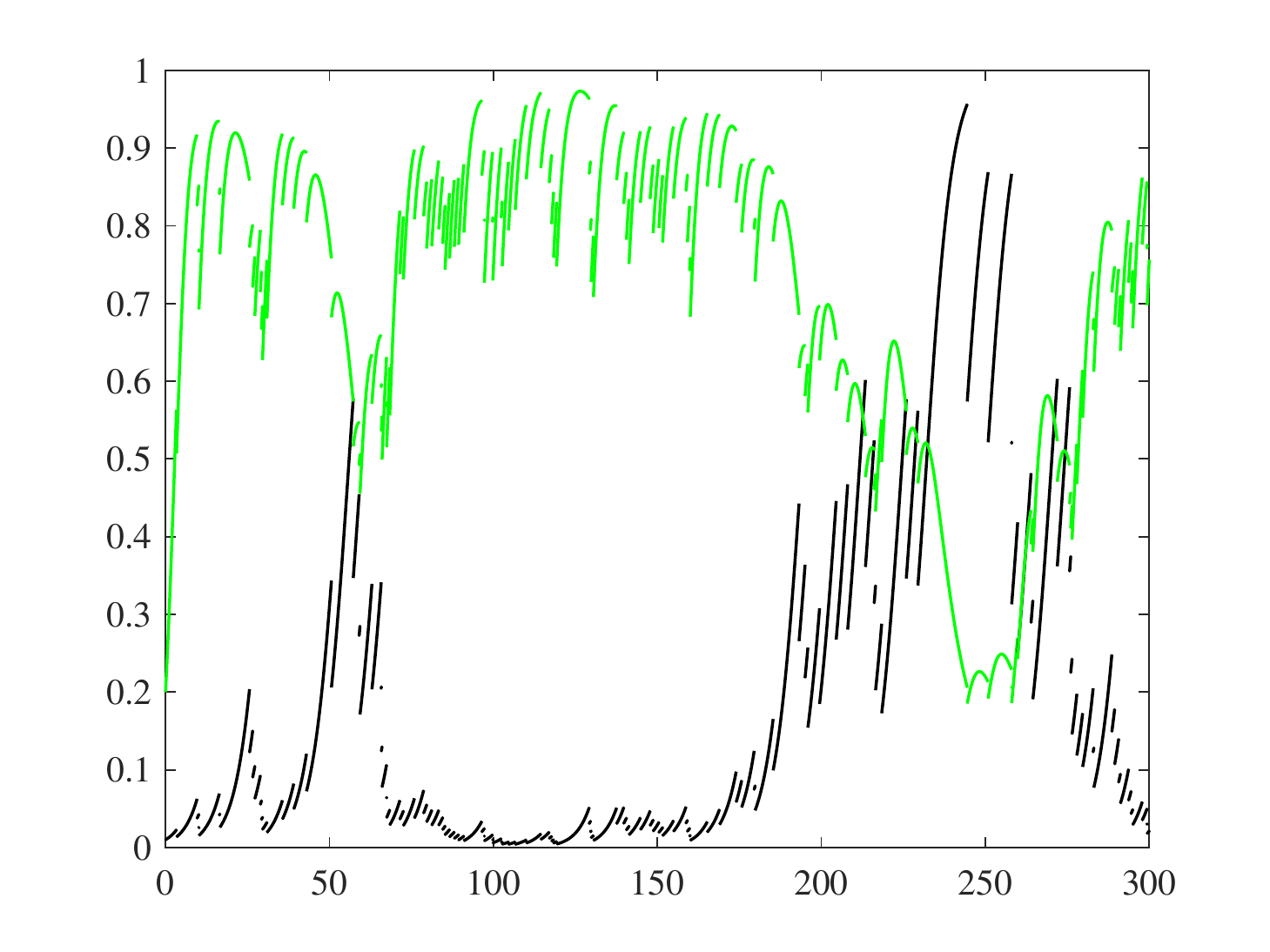}
 \caption{Sample trajectories of the stochastic process in \eqref{basic.eq}--\eqref{loss.eq} with parameter values $r_w=0.25$, $r_g=0.5$, $M_w=0.4$, $M_g=0.1$, $\lambda(w,g)=g$ and initial condition $w_0=0.01$, $g_0=0.2$
 } \label{wg.fire}
\end{figure}

Instead of using actual amount of biomasses we will relate in our model to ratios of these amounts to maximal capacities of wood and grass, respectively:
\[
w(t)=\frac{W(t)}{K_W},\quad g(t)=\frac{G(t)}{K_G}.
\]
Thus $w(t)$ and $g(t)$ take values in the unit interval, i.e. $0\leqslant w(t), g(t)\leqslant1$ for any time $t$.
 We now allow disturbances of the growth of biomasses due to fires occurring at random times $(t_n)_{n\geqslant 1}$. Let $t_0=0$ and $w(t_0)=w_0$, $g(t_0)=g_0$, where $w_0,g_0\in [0,1]$ are arbitrary. In periods
between fire occurrences the growth of normalized tree and grass biomasses is modeled with
\begin{equation}\label{basic.eq}
\left\{\begin{array}{lll}
w'(t)=r_w  w(t)\big(1-w(t)\big),\\
g'(t)=r_g  g(t)\big(1-g(t)-w(t)\big),
\end{array} \right.
\end{equation}
for $t\in(t_{n}, t_{n+1})$, $n\ge 0$, and a sequence of random variables $(\tau_n)$ such that
\begin{equation}\label{tau}
\left\{\begin{array}{lll}
t_{n+1}=t_{n}+\tau_{n+1},\\
\Pr\big(\tau_{n+1}>t|w(t_{n})=w_{n}, g(t_{n})=g_{n}\big)=e^{-\int_0^t \lambda(\pi_s(w_n,g_n))ds},
\end{array} \right.
\end{equation}
where $\pi_t(w_n,g_n)=(w(t),g(t))$ is the solution of  (\ref{basic.eq}) with initial condition $(w_n,g_n)$ and $\lambda$ is a nonnegative bounded continuous function.
At each time $t_{n+1}$ the loss of biomasses is given by
\begin{equation}\label{loss.eq}
\left\{\begin{array}{lll}
w(t_{n+1})=w(t_{n+1}^-)-M_w \, w(t_{n+1}^-),\\
g(t_{n+1})=g(t_{n+1}^-)-M_g \, g(t_{n+1}^-),\quad n\ge 0,
\end{array} \right.
\end{equation}
where $ M_w,M_g\in (0,1)$ are constants, $v(t^-)=\lim_{s\rightarrow t^-}v(s)$ for $v\in\{w, g\}$.  We assume that the function $\lambda\colon [0,1]^2\to \mathbb{R}_+$ satisfies
\begin{equation}\label{lambda}
\lambda(w,0)=0, \quad w\ge 0,\quad \lambda(w,g)>0\quad \text{for } w\geqslant 0,\ g>0.
\end{equation}
In Figure \ref{fire-no-fire} we display  graphs of wood and grass biomasses in time, without and with fires. A sample behavior of the overall system in the long run including
losses due to random fires  is shown in Figure \ref{wg.fire}.

The process $\xi(t)=(w(t),g(t))$, $t\ge 0$,  with $w,g$ as in \eqref{basic.eq}--\eqref{loss.eq},  is a piecewise deterministic Markov process (\cite{davis84}) with state space $[0,1]^2$. It is an example of a flow with jumps as presented in~\cite[Section 4.2.4]{rudnickityran17}.
We describe the jumps of the stochastic process by a linear transformation $S$ mapping $(w,g)\longmapsto S(w,g)$, where
\begin{equation}\label{d:S}
S(w,g)=((1-M_w)w,(1-M_g)g),\quad (w,g)\in [0,1]^2.
\end{equation}
Let $p(t,w,g)$ be the probability density of $(w(t),g(t))$, i.e. $p$ is nonnegative, Borel measurable, and satisfies
\[
\Pr((w(t),g(t))\in B)=\int_B p(t,w,g)dwdg
\]
for any Borel subset of $[0,1]^2$ with the integral being equal to one for $B=[0,1]^2$. Then $p$ is a solution of the following Fokker-Planck type equation \begin{multline}\label{fokk.pla.eq}
\frac{\partial p(t,w,g)}{\partial t}+\frac{\partial\big(r_ww(1-w)p(t,w,g)\big)}{\partial w}+\frac{\partial\big(r_gg(1-g-w)p(t,w,g)\big)}{\partial g}\\=-\lambda(w,g)p(t,w,g)
+\frac{\lambda\big(S^{-1}(w,g))\big)p\big(t, S^{-1}(w,g)\big)}{(1-M_w)(1-M_g)},
\end{multline}
where $S^{-1}$ is the inverse of the transformation $S$ defined in \eqref{d:S}. Equation \eqref{fokk.pla.eq} is supplemented
with initial condition
\begin{equation}\label{fokk.pla.eq.ic}
p(0,w,g)=f(w,g),\quad \text{where } \int_{0}^1\int_0^1 f(w,g)\,dwdg=1
\end{equation}
and $f$ is a nonnegative Borel measurable function, so that $f$ is the probability density of $(w(0),g(0))$.
We have the following result - its proof will be given in the next section.

\begin{thm}\label{main}
There exists a unique density $p_*(w,g)$ which is a stationary solution of \eqref{fokk.pla.eq}. Moreover, every solution of \eqref{fokk.pla.eq}--\eqref{fokk.pla.eq.ic} converges to $p_*$, i.e.
\[
\lim_{t\rightarrow\infty}\int_0^1\int_0^1 \big|p(t,w,g)-p_*(w,g)\big|dwdg=0.
\]
\end{thm}
\begin{remark}\label{r:boun} Let $w_0,g_0\in [0,1]$ and $w(t),g(t)$ be as in \eqref{basic.eq}--\eqref{loss.eq}.
If $w_0>0$ and $g_0=0$ then $g(t)=0$ for all $t>0$. In this case, assumption \eqref{lambda} implies that fire can not occur when there is no grass biomass.
Hence, $w$ is defined for all $t$  as the solution of the differential equation $w'(t)=r_w w(t)(1-w(t))$ with initial condition $w(0)=w_0$. Thus  $w(t)>0$ for all $t>0$ and $w(t)$ converges to $1$ as $t\to \infty$. Consequently, the point measure $\delta_{\{(1,0)\}}$ is an invariant measure for the process $\xi$.
Similarly, if $w_0=0$ and $g_0=0$ then $w(t)=0$ and $g(t)=0$ for all $t> 0$. Thus also the point measure $\delta_{\{(0,0)\}}$ is an invariant measure for the process $\xi$. Finally, if $w_0=0$ and $g_0>0$ then $w(t)=0$ and $g(t)>0$ for all $t\ge 0$. In this case, the process $\xi$ has an invariant distribution which is a product of $\delta_{\{0\}}$ and an absolutely continuous measure, see Remark~\ref{r:inv}.
\end{remark}

\begin{remark}
If instead of \eqref{tau} we have  $t_{n+1}=t_n+\tau$, $n\ge 0$ where $\tau$ is a constant then such a model is an example of an impulsive  system \cite{tamen2016tree,yatat2017impulsive}.
\end{remark}

\section{Existence and uniqueness of tree and grass biomasses distribution}
Methods in this section are mostly taken from the book \cite{rudnickityran17}.
To prove Theorem \ref{main} we use the method from
\cite[Section 6.3.2]{rudnickityran17}. We begin by  recalling some notions for stochastic semigroups.
Let the triple $(X,\Sigma,m)$ be a $\sigma$-finite measure space.
Denote by $D$ the subset of the space
$L^1=L^1(X,\Sigma,m)$ which contains all
densities
\[
D=\{f\in  L^1:  f\ge 0,\; \|f\|=1\}.
\]
A linear mapping $P\colon L^1\to L^1$ is called a
\emph{Markov} or \emph{stochastic operator\,} if $P(D)\subset D$. A family $\{P(t)\}_{t\ge0}$ of stochastic operators
which satisfies conditions:
\begin{enumerate}
\item  $P(0)=\operatorname{id}$, $P(t+s)=P(t) P(s)$ for
$s,\,t\ge0$,
\item  for each $f\in L^1$ the function
$t\mapsto P(t)f$ is continuous,
\end{enumerate}
 is called a
\emph{stochastic semigroup}.

Consider a stochastic semigroup
$\{P(t)\}_{t\ge0}$.
A density $f_*$
is called {\it invariant\,} if $P(t)f_*=f_*$ for each $t>0$.
The stochastic semigroup
$\{P(t)\}_{t\ge 0}$
is called {\it asymptotically
stable\,} if there is an invariant density
$f_*$ such that
\[
\lim _{t\to\infty}\|P(t)f-f_*\|=0 \quad \text{for}\quad
f\in D.
\]

We will use a result of Pichór and Rudnicki \cite{pichorrudnicki18}  (see also \cite[Theorem 5.6]{rudnickityran17}) which requires the following conditions:
\begin{enumerate}
\item[(K)] For every $y_0\in X$ there exist $\varepsilon >0$, $t>0$,
and a measurable function
$\eta\ge 0$ such that $\int \eta(x)\, m(dx)>0$ and
\begin{equation*}
\label{w-eta2}
P(t)f(x)\ge \eta(x)\int_{B(y_0,\varepsilon)} f(y)\,m (dy),
\end{equation*}
where
$B(y_0,\varepsilon)=\{y\in X:\,\,\rho(y,y_0)<\varepsilon\}$.
\item[(WI)]  There exists a point $x_0\in X$ such that for each $\varepsilon >0$ and for each density $f$ we have
\begin{equation*}
\label{k:WI}
\int\limits_{B(x_0,\varepsilon)} P(t)f(x)\,m(dx)>0\quad\textrm{for some $t=t(\varepsilon,f)>0$}.
\end{equation*}
\item[(WT)]  There exists $\kappa>0$ such that
\begin{equation*}
\label{k:T}
\sup\limits_{F\in \mathcal F}\limsup_{t\to\infty} \int_F P(t)f(x)\,m(dx)\ge \kappa
\end{equation*}
for $f\in D_0$, where $D_0$ is a dense subset of $D$ and $\mathcal F$ is the family of all compact subsets of $X$.
\end{enumerate}

\begin{thm}
\label{cor:st}
Let $\{P(t)\}_{t\ge0}$ be a stochastic semigroup on
$L^1(X,\Sigma,m)$, where
$X$ is a separable metric space,  $\Sigma$ is the $\sigma$-algebra of Borel subsets of $X$, and $m$ is a $\sigma$-finite measure.
Assume that $\{P(t)\}_{t\ge0}$ satisfies conditions {\rm (K)}, {\rm (WI)}, and {\rm (WT)}.
Then the semigroup $\{P(t)\}_{t\ge0}$ is asymptotically stable.
\end{thm}

Now we are ready to prove the main theorem.

\begin{proof}[of Theorem \ref{main}]
Let $X=(0,1]^2$  and $m$ be the two-dimensional Lebesgue measure on $X$. It follows from \cite[Section 4.2.4]{rudnickityran17} that the process $\xi(t)$, $t\ge 0$, induces a stochastic semigroup $\{P(t)\}_{t\ge 0}$  on $L^1=L^1(X,\Sigma,m)$ and that the solution of \eqref{fokk.pla.eq}--\eqref{fokk.pla.eq.ic} is given by $p(t,w,g)=P(t)f(w,g)$, $t\ge 0$, $(w,g)\in X$.  To apply Theorem~\ref{cor:st} we need to check conditions (K), (WI), and (WT).

We first show that condition (WT) holds.
The extended generator  $\widetilde{L}$ of the process $\xi$ is of the form
\[
\widetilde{L}V(x)=\langle b(x),\mathrm{grad}V(x)\rangle +\varphi(x)(V(S(x))-V(x))\quad \text{for }x=(w,g),
\]
where $\mathrm{grad}V(x)$ is the gradient of $V(x)$ and $b(x)$ is the vector with coordinates
\[
b_1(x)=r_w w(1-w),\quad b_2(x)=r_gg(1-w-g), \quad x=(w,g).
\]
The domain $\mathcal{D}(\widetilde{L})$  of the extended generator $\widetilde{L}$ (see \cite{davis84} or \cite[Section 2.3.6]{rudnickityran17})  contains the set of functions $V\colon X\to \mathbb{R}$ such that
for each $x\in X$ the function $t\mapsto V(\pi_t(x))$ is absolutely continuous and  for each $t\ge 0$, $x\in X$, we have \[
\mathbb{E}\Big(\sum_{t_n\le t}\big|V(\xi(t_n))-V(\xi(t_n^{-}))\big|\Big|\xi(0)=x\Big)<\infty.
\]
Let $V(w,g)=-\log w -\log g$. Since we have
$V(\xi(t_n))-V(\xi(t_n^{-}))=-\log(1-M_w)-\log(1-M_g)$ for any $n$, we see that $V$ belongs to $\mathcal{D}(\widetilde{L})$  and that
\[
\widetilde{L}V(w,g)=-r_w (1-w)-r_g(1-w-g)-\lambda(w,g)(\log(1-M_w)+\log(1-M_g)).
\]
The function $\widetilde{L}V$ is bounded on $(0,1]^2$ and $\widetilde{L}V(w,g)\to -r_w-r_g$ as $\|(w,g)\|\to 0$, where $\|\cdot\|$ denotes a norm in $\mathbb{R}^2$. Thus we can find a $\delta\in (0,1)$  such that $\widetilde{L}V(w,g)\le -(r_w+r_g)/2$ for $\|(w,g)\|< \delta$.
Moreover, we have
\[
\int_X \widetilde{L}V(x)f(x)m(dx)=\int_X V(x)Af(x)m(dx),\quad f\in \mathcal{D}(A)\cap \mathcal{D}_V,
\]
where $f\in D_V$ iff $\int_{X}V(x)|f(x)|m(dx)<\infty$ and $(A,\mathcal{D}(A))$ is the generator of the semigroup $\{P(t)\}_{t\ge 0}$.
We conclude that $V$ is a Hasminski\u{\i} function for the semigroup $\{P(t)\}_{t\ge 0}$ and the  compact set $F=\{(w,g)\in X: \|(w,g)\|\ge  \delta\}$ implying that condition (WT) holds, by \cite[Corollary 5.8]{rudnickityran17}.

To check condition (K) take  $x_0\in X$ and define
 $x_1=S(x_{0})$, $x_2=S(x_1)$, 
\begin{equation}
\label{a-cs-1}
v_1=S'(x_{1})S'(x_{0})b(x_{0})
-b(x_{2}),\quad v_2=S'(x_{1})b(x_{1})
-b(x_{2}).
\end{equation}
Since $S$ is a linear transformation, we have
\[
v_1=S^2(b(x_{0}))
-b(S^2(x_{0})),\quad v_2=S(b(S(x_{0}))
-b(S^2(x_{0})),
\]
where $S^2(x)=S(S(x))$.
It is easily seen that vectors $v_1$ and $v_2$ are linearly independent for each $x_0\in X$. Since the function $\lambda$ is strictly positive on $(0,1]^2$, we conclude  that condition (K) holds (see e.g. \cite[Section 6.3.2]{rudnickityran17} or \cite[Section 4]{biedrzyckatyran}).

Observe that condition (WI) holds once we show that there exists $x_0$ such that for each $\varepsilon>0$ and $x\in X$ we can find $n$ and times $s_1,\ldots,s_n,s_{n+1}>0$ such that  $\pi_{s_{n+1}}(x_n)\in B(x_0,\varepsilon)$,
where
\begin{equation}\label{d:xn}
x_{n}=S(\pi_{s_{n}}(\ldots S(\pi_{s_1}x))).
\end{equation}
The point $(0,1)$ is a saddle point for the two-dimensional system \eqref{basic.eq} considered on $\mathbb{R}^2$. Its stable manifold is the set $\{(0,g): g>0\}$ and its unstable manifold contains a curve joining the point $(0,1)$ with the stable point $(1,0)$, see Figure~\ref{fig-phase}. Let us take  $x_0\in X$ from this curve lying close to the point $(1,0)$. For any point $x\in X$ we can find $n$ and $s_1,\ldots,s_n>0$  such that $x_{n}$ defined as in \eqref{d:xn} is as close to  $(0,0)$ as is needed. Since $\pi_s y\to (1,0)$ for $y\in (0,1)^2$, we can find $s_{n+1}$ such that $\pi_{s_{n+1}} x_n\in B(x_0,\varepsilon)$, which completes the proof.
\end{proof}

\begin{remark}\label{r:inv}
The process $\xi$ restricted to the set $\{(0,g):g\in(0,1]\}$, considered with measure being the product of $\delta_{\{0\}}$ and the Lebesgue, induces a stochastic semigroup on $L^1(\{0\}\times(0,1])$. Using the same type of argument as in the proof of Theorem \ref{main} it can be shown that this semigroup  satisfies conditions (K), (WI) and (WT), thus this semigroup is asymptotically stable, implying the existence of the invariant measure mentioned in Remark~\ref{r:boun}.
\end{remark}

\section{Discussion}
We showed that there exists unique, absolutely continuous with respect to the two-dimensional Lebesgue measure, stationary distribution for positive amount of grass and wood biomasses. The stationary density is strictly positive in the region  bounded by the axes and the unstable manifold of the point $(0,1)$, in particular in a neighbourhood of the line $\{(w,1-w):w\in (0,1]\}$, showing that the coexistence of trees and grass is possible.
Finding the actual shape of this distribution, numerical analysis, and further improvements of the model by adding more coefficients reflecting real-world factors regulating savanna biomasses we leave for future work.

Moreover such analysis can be implemented in models describing different phenomena involving random fires, such as impact of forest fires on population of pines and bark beetles. Modeling attempts usually are deterministic (see e.g. \cite{chen2014model}) and hence could benefit from involving stochastic nature of fire.

\end{document}